\documentclass{amsart}
\usepackage{amssymb}

\numberwithin{equation}{section}

\newtheorem{propo}{Proposition}[section]
\newtheorem{corol}[propo]{Corollary}
\newtheorem{theor}[propo]{Theorem}
\newtheorem{lemma}[propo]{Lemma}
\theoremstyle{definition}
\newtheorem{defin}[propo]{Definition}
\newtheorem{examp}[propo]{Example}
\theoremstyle{remark}
\newtheorem{remar}[propo]{Remark}

\newcommand{\ord}{\operatorname{ord}}
\newcommand{\Hom}{\operatorname{Hom}}
\newcommand{\bfg}{\mathbf{g}}
\newcommand{\bfch}{\mathbf{\chi}}
\newcommand{\bflm}{\mathbf{\lambda}}
\newcommand{\cou}{\varepsilon}
\newcommand{\lmb}{\lambda}
\newcommand{\Ss}{\mathcal{S}}
\newcommand{\wt}{\widetilde}

\newcommand{\bZ}{\mathbb{Z}}
\newcommand{\bN}{\mathbb{N}}
\newcommand{\bk}{\mathbf{k}}
\newcommand{\al}{\alpha}
\newcommand{\com}[2]{#1_{(#2)}}                 
\newcommand{\so}{\Rightarrow}

\newcommand{\id}{\operatorname{id}}
\begin{document}
\title[Decomposition of some pointed Hopf algebras]{Decomposition of some pointed Hopf algebras given by the canonical Nakayama automorphism}
\author{M. Gra\~na}
\author{J.A. Guccione}
\author{J.J. Guccione}
\address{%
M.G., J.A.G., J.J.G.:\newline\indent
    Depto de Matem\'atica - FCEyN\newline\indent
    Universidad de Buenos Aires\newline\indent
    Pab. I - Ciudad Universitaria\newline\indent
    1428 - Buenos Aires - Argentina\newline
	\email{M.G.}{matiasg@dm.uba.ar}
	\email{J.A.G.}{vander@dm.uba.ar}
	\email{J.J.G.}{jjgucci@dm.uba.ar}}
\thanks{This work was partially supported by CONICET, PICT-02 12330, UBA X294}
\begin{abstract}
Every finite dimensional Hopf algebra is a Frobenius algebra, with Frobenius homomorphism
given by an integral. The Nakayama automorphism determined by it yields a decomposition with
degrees in a cyclic group. For a family of pointed Hopf algebras, we determine necessary and
sufficient conditions for this decomposition to be strongly graded.
\end{abstract}


\maketitle

\section{Introduction}
Let $\bk$ be a field, $A$ a finite dimensional $\bk$-algebra and $DA$ the dual space
$\Hom_\bk(A,\bk)$, endowed with the usual $A$-bimodule structure. Recall that $A$ is said to
be a Frobenius algebra if there exists a linear form $\varphi\colon A\to \bk$, such that the
map $A\to DA$, defined by $x\mapsto x\varphi$, is a left $A$-module isomorphism. This linear
form $\varphi\colon A\to \bk$ is called a Frobenius homomorphism. It is well known that this
is equivalent to say that the map $x\mapsto \varphi x$, from $A$ to $DA$, is an isomorphism of
right $A$-modules. From this it follows easily that there exists an automorphism $\rho$ of
$A$, called the Nakayama automorphism of $A$ with respect to $\varphi$, such that $x\varphi =
\varphi \rho(x)$, for all $x\in A$. It is easy to check that a linear form $\wt{\varphi}\colon
A\to \bk$ is another Frobenius homomorphism if and only if there exists an invertible element
$x$ in $A$, such that $\wt{\varphi} = x\varphi$. It is also easy to check that the Nakayama
automorphism of $A$ with respect to $\wt{\varphi}$ is the map given by $a\mapsto
\rho(x)^{-1}\rho(a)\rho(x)$.

\smallskip

Let $A$ be a Frobenius $\bk$-algebra, $\varphi\colon A\to \bk$ a Frobenius homomorphism and
$\rho:A\to A$ the Nakayama automorphism of $A$ with respect to $\varphi$.

\begin{defin}
	We say that $\rho$ has order $m\in\bN$ and we write $\ord_{\rho} = m$,
	if $\rho^m = id_A$ and $\rho^r\ne id_A$, for all $r<m$.
\end{defin}

Assume that $\rho$ has finite order and that $\bk$ has a primitive $\ord_{\rho}$-th root of
unity $\omega$. For $n\in\bN$, let $C_n$ be the group of $n$-th roots of unity in $\bk$. Since
the polynomial $X^{\ord_{\rho}}-1$ has distinct roots $\omega^i$ ($0\le i<\ord_{\rho}$), the
algebra $A$ becomes a $C_{\ord_{\rho}}$-graded algebra
\begin{equation}\label{eq:lad}
A = A_{\omega^0}\oplus \cdots\oplus A_{\omega^{\ord_{\rho}-1}},
	\quad\text{where }A_z = \{a\in A\;:\;\rho(a) = za\}.
\end{equation}
As it is well known, every finite dimensional Hopf algebra $H$ is Frobenius, being a Frobenius
homomorphism any nonzero right integral $\varphi\in H^*$. Let $t$ be a nonzero right integral of
$H$. By \cite[Proposition~3.6]{S}, the compositional inverse of the Nakayama map $\rho$ with
respect to $\varphi$, is given by
$$
\rho^{-1}(h) = \al(h_{(1)})\Ss^{-2}(h_{(2)}),
$$
where $\al\in H^*$ is the modular element of $H$, defined by $at = \al(a)t$ (note that the
automorphism of Nakayama considered in \cite{S} is the compositional inverse of the one
considered by us). Using this formula and that $\al\circ \Ss^2 = \al$, it is easy to check
that $\rho(h) = \al(\Ss(h_{(1)}))\Ss^2(h_{(2)})$, and more generally, that
\begin{equation}\label{eq:rhol}
	\rho^l(h) = \al^{*l}(\Ss(h_{(1)}))\Ss^{2l}(h_{(2)}),
\end{equation}
where $\al^{*l}$ denotes the $l$-fold convolution product of $\al$. Since $\al$ has finite
order with respect to the convolution product and, by the Radford formula for $\Ss^4$ (see
\cite{R} or \cite[Theorem~3.8]{S}), the antipode $\Ss$ has finite order with respect to
composition, the automorphism $\rho$ has finite order. So, finite dimensional Hopf algebras
are examples of the situation considered above.

Notice that by \eqref{eq:rhol}, if $\rho^l = \id$, then $\alpha^{*l} = \cou$ and then
$\Ss^{2l} = \id$. The converse is obvious. So, the order of $\rho$ is the lcm between those of
$\alpha$ and $\Ss^2$. In particular, the number of terms in the decomposition associated with
$\Ss^2$ divides that in the one associated with $\rho$. Also, from \eqref{eq:rhol} we get that
$\rho = \Ss^2$ if and only if $H$ is unimodular.

\smallskip

The main aim of the present work is to determine conditions for decomposition~\eqref{eq:lad}
to be strongly graded. Besides the fact that the theory for algebras which are strongly graded
over a group is well developed (see for instance \cite{galgs}), our interest on this problem
originally came from the homological results in \cite{gg}.

The decomposition using $\Ss^2$ instead of $\rho$ was considered in \cite{rs}. We show below that if
$\Ss^2\neq\id$, then this decomposition is not strongly graded. On the other hand, as shown in \cite{rs},
under suitable assumptions its homogeneous components are equidimensional.  It is an interesting problem to
know whether a similar thing happens with the decomposition associated with $\rho$. For instance, all the
liftings of Quantum Linear Spaces have equidimensional decompositions, as shown in Remark~\ref{rm:eqc}.

\section{The unimodular case}
Let $H$ be a finite dimensional Hopf algebra with antipode $\Ss$. In this brief section we first show that the
decomposition of $H$ associated with $\Ss^2$ is not strongly graded, unless $\Ss^2=\id$ (this applies in
particular to decomposition~\eqref{eq:lad} when $H$ is unimodular and $\ord_{\rho}>1$). We finish by giving a
characterization of unimodular Hopf algebras in terms of decomposition~\eqref{eq:lad}.

\begin{lemma}\label{lm:nsg}
Let $H$ be a finite dimensional Hopf algebra. Suppose $H=\bigoplus_{g\in G}H_g$ is a graduation over a group.
Assume there exists $g\in G$ such that $\varepsilon(H_g)=0$. Then the decomposition is not strongly graded.
\end{lemma}
\begin{proof}
	Suppose the decomposition is strongly graded. Then there are elements $a_i\in H_g$ and $b_i\in H_{g^{-1}}$
	such that $1=\sum_ia_ib_i$. Then $1=\varepsilon(1)=\sum_i\varepsilon(a_i)\varepsilon(b_i)$, a
	contradiction.
\end{proof}

\begin{corol}
	Assume that $\Ss^2\neq\id$ and that
	$$
		H=\bigoplus_{z\in\bk^*}H_z, \quad\text{where } H_z=\{h\in H\;:\;\Ss^2(h)=zh\}.
	$$
	Then this decomposition is not strongly graded.
\end{corol}
\begin{proof}
	Since $\varepsilon\circ\Ss^2=\varepsilon$, then $\varepsilon(H_z)=0$ for all $z\neq 1$.
\end{proof}

Let now $\varphi\in \int_{H^*}^r$ and $\Gamma\in \int_H^l$, such that $\langle \varphi,\Gamma \rangle = 1$,
and let $\alpha\colon H\to \bk$ be the modular map associated with $t = S(\Gamma)$. Let $\rho$ be the Nakayama
automorphism associated with $\varphi$. Assume that $\bk$ has a root of unity $\omega$ of order $\ord_\rho$.
We consider the decomposition associated with $\rho$, as in \eqref{eq:lad}
\begin{equation}\label{umc1}
	H = H_{\omega^0}\oplus\cdots\oplus H_{\omega^{\ord_{\rho} - 1}}.
\end{equation}

\begin{corol}
	If $H$ is unimodular and $\Ss^2\neq\id$, then the decomposition~\eqref{umc1} is not strongly graded.
	\hfill\qed
\end{corol}

\begin{propo}
If $h\in H_{\omega^i}$, then $\alpha(h) = \omega^{-i}\epsilon(h)$.
\end{propo}

\begin{proof}
In the proof of \cite[Proposition~3.6]{S} it was shown that $\langle\varphi, t \rangle = 1$.
Then
$$
\epsilon(h) = \epsilon(h)\langle \varphi, t \rangle = \langle \varphi, th \rangle = \langle
\varphi, \rho(h)t \rangle = \langle\varphi, \omega^i ht \rangle = \omega^i\alpha(h).
$$
So, $\alpha(h) = \omega^{-i}\epsilon(h)$, as we want.
\end{proof}

\begin{corol}\label{co:huiff}
$H$ is unimodular if and only if $H_{\omega^i}\subseteq \ker(\epsilon)$, for all $i>0$
\end{corol}

\begin{proof}
$\Rightarrow$):\enspace For $h\in H_{\omega^i}$, we have $\epsilon(h) = \alpha(h) =
\omega^{-i} \epsilon(h)$ and so $\epsilon(h) = 0$, since $\omega^i\neq 1$.

\smallskip

\noindent $\Leftarrow$):\enspace For $h\in H_{\omega^i}$ with $i>0$, we have $\alpha(h) =
\omega^{-i} \epsilon(h) = 0 = \epsilon(h)$ and, for $h\in H_{\omega^0}$, we also have
$\alpha(h) = \omega^0 \epsilon(h) = \epsilon(h)$.
\end{proof}

\section{Bosonizations of Nichols algebras of diagonal type}
Let $G$ be a finite abelian group, $\bfg = g_1,\ldots,g_n\in G$ a sequence of elements in $G$ and
$\bfch = \chi_1,\ldots,\chi_n\in\hat G$ a sequence of characters of $G$. Set $q_{ij} = \chi_j(g_i)$. Let $V$
be the vector space with basis $\{x_1,\ldots,x_n\}$ and let $R = \mathfrak{B}(V)$ be the Nichols algebra
generated by $(V,c)$, where $c$ is the braiding given by $c(x_i\otimes x_j) = q_{ij}x_j\otimes x_i$. Let
$T_c(V)$ be the tensor algebra generated by $V$, endowed with the unique braided Hopf algebra structure such
that the elements $x_i$ are primitive and whose braiding extends $c$. Then, $R$ is the quotient of $T_c(V)$ by
the ideal generated by those primitive homogeneous elements with degree $\ge 2$. See \cite{ag,as} for the
definition and main properties of Nichols algebras. We give here one of the possible equivalent definitions:
assume that $R$ is finite-dimensional and let $t_0\in R$ be a nonzero homogeneous element of greatest degree.
Let $H = H(\bfg,\bfch) = R\#\bk G$ be the bosonization of $R$ (this is an alternative presentation for the
algebras considered by Nichols in \cite{n}). We have:
\begin{align*}
	& \Delta(g_i) = g_i\otimes g_i, && \Delta(x_i) = g_i\otimes x_i + x_i\otimes 1,\\
	&\Ss(g) = g^{-1},               && \Ss(x_i) = -g_i^{-1}x_i,\\
	& \Ss^2(g) = g,                 && \Ss^2(x_i) = g_i^{-1}x_ig_i = q_{ii}^{-1}x_i.
\end{align*}
The element $t_0\sum_{g\in G}g$ is a non zero right integral in $H$. Let $\alpha$ be the
modular element associated with it. Thence, $\alpha(x_i) = 0$ for all $i$,
and $\alpha_{\vert G}$ is determined by $gt_0 = \alpha(g)t_0g$. Thus,
$$
\rho(g) = \alpha(g^{-1}) g \quad\text{and}\quad \rho(x_i) = \alpha(g_i^{-1}) q_{ii}^{-1} x_i.
$$
This implies that the nonzero monomials $x_{i_1}\cdots x_{i_\ell}g$ are a set of eigenvectors
for $\rho$ (which generate $H$ as a $k$-vector space). Consider the subgroups
$$
L_1 = \langle q_{11},\ldots,q_{nn},\alpha(G)\rangle\quad\text{and}\quad L_2 = \alpha(G),
$$
of $k^*$. Since $\rho(x_ig_i^{-1}) = q_{ii}^{-1}x_ig_i^{-1}$ and $\rho(g) = \alpha(g)$, the
group $L_1$ is the set of eigenvalues of $\rho$. Consequently $\bk$ has a primitive $\ord_{\rho}$-th root of
unity. As in the introduction, we decompose
\[
H = \bigoplus_{\omega\in L_1} H_{\omega},\quad\text{where } H_{\omega} = \{h\in H\;:\;\rho(h) =
\omega h\}.
\]

\begin{propo}\label{pr:gralc1} The following are equivalent:
\begin{enumerate}
\item\label{pr:gralit1} $\bigoplus_{\omega\in L_1} H_{\omega}$ is strongly graded,
\item\label{pr:gralit2} $L_1 = L_2$,
\item\label{pr:gralit3} Each $H_{\omega}$ contains an element in $G$.
\item\label{pr:gralit4} $H$ is a crossed product $H_1\ltimes \bk L_1$.
\end{enumerate}
\end{propo}

\begin{proof} It is clear that
\eqref{pr:gralit2}$\so$\eqref{pr:gralit3}$\so$\eqref{pr:gralit4}$\so$\eqref{pr:gralit1}. We
now prove that \eqref{pr:gralit1}$\so$\eqref{pr:gralit2}. Notice that $H$ is also
$\bN_0$-graded by $\deg(g) = 0$ for all $g\in G$, and $\deg(x_i) = 1$.
Call $H = \bigoplus_{i\in\bN}H^i$ this decomposition. Since each
$H_{\omega}$ is spanned by elements which are homogeneous with respect to the previous
decomposition, we have:
$$
H = \bigoplus_{i\ge 0,\;\omega\in L_1} H_{\omega}^i,\quad\text{where } H_{\omega}^i =
H_{\omega}\cap H^i.
$$
So, if $\bigoplus_{\omega\in L_1}H_{\omega}$ is strongly graded, then each $H_{\omega}$ must
contain nonzero elements in $H^0$. Since $H^0\subseteq\bigoplus_{\omega\in L_2} H_{\omega}$,
we must have $L_1 = L_2$.
\end{proof}

\subsection*{Quantum Linear Spaces}
If the sequence of characters $\chi$ satisfies
\begin{itemize}
\item $\chi_i(g_i)\neq 1$,
\item $\chi_i(g_j)\chi_j(g_i) = 1$ for $i\neq j$,
\end{itemize}
then $H(\bfg,\bfch)$ is the Quantum Linear Space with generators $G$ and $x_1,\ldots,x_n$,
subject to the following relations:
\begin{itemize}
\item $gx_i = \chi_i(g)x_ig$,
\item $x_ix_j = q_{ij}x_jx_i$,
\item $x_i^{m_i} = 0$,
\end{itemize}
where $m_i = \ord(q_{ii})$.  For these sort of algebras it is possible to give an explicit formula for $\rho$.
In fact, the element $t = x_1^{m_1-1}\cdots x_n^{m_n-1}\sum_{g\in G}g$ is a right integral in $H$. Using this
integral, it is easy to check that $\alpha(g) = \chi_1^{m_1-1}(g) \cdots\chi_n^{m_n-1}(g)$.
In particular, $\alpha(g_i) = q_{i1}^{m_1-1}\cdots q_{in}^{m_n-1}$.
A straightforward computation, using that $\rho(g) = \alpha(g^{-1})g$ and
$\rho(x_i) = \alpha(g_i^{-1})q_{ii}^{-1}x_i$, shows that
\[
\rho(x_1^{r_1}\cdots x_n^{r_n}g) = \prod_{1\le i<j\le n} q_{ij}^{(1-m_j)r_i-(1-m_i)r_j}
\alpha(g^{-1}) x_1^{r_1}\cdots x_n^{r_n}g.
\]
Proposition~\ref{pr:gralc1} applies to this family of algebras.

\begin{examp} Let $\bk$ be a field of characteristic $\neq 2$ and let $G = \{1,g\}$. Set $g_i=g$
	and $\chi_i(g) = -1$ for $i\in \{1,\ldots,n\}$. Then, $q_{ij} = -1$ for all $i,j$,
	and $\alpha(g) = (-1)^n$. In this case, the algebra $H$ is generated by $g,x_1,\ldots,x_n$
	subject to relations
	\begin{itemize}
	\item $g^2 = 1$,
	\item $x_i^2 = 0$,
	\item $x_ix_j = -x_jx_i$,
	\item $gx_i = -x_ig$.
	\end{itemize}
	By Proposition~\ref{pr:gralc1}, we know that $H$ is strongly graded if and only if $n$ is odd.
\end{examp}

\section{Liftings of Quantum Linear Spaces}
In this section we consider a generalization of Quantum Linear Spaces: that of their liftings.
As above, $G$ is a finite abelian group, $\bfg = g_1,\ldots,g_n\in G$ is a sequence of
elements in $G$ and $\bfch = \chi_1,\ldots,\chi_n\in\hat G$ is a sequence of characters of
$G$, such that
\begin{eqnarray}
&&\chi_i(g_i)\neq 1, \\
&&\chi_i(g_j)\chi_j(g_i) = 1, \text{ for } i\neq j.
\end{eqnarray}
Again, let $q_{ij} = \chi_j(g_i)$ and let $m_i = \ord(q_{ii})$. Let now $\lmb_i\in\bk$
and $\lmb_{ij}\in\bk$ for $i\neq j$ be such that
$$
\lmb_i(\chi_i^{m_i}-\cou) = \lmb_{ij}(\chi_i\chi_j-\cou) = 0.
$$
Suppose that $\lmb_{ij}+q_{ij}\lmb_{ji} = 0$ whenever $i\neq j$. The \emph{lifting} of the
quantum affine space associated with this data is the algebra $H = H(\bfg,\bfch,\bflm)$, with
generators $G$ and $x_1,\ldots,x_n$, subject to the following relations:
\begin{eqnarray}
&&gx_i = \chi_i(g)x_ig,\\
&&x_ix_j = q_{ij}x_jx_i+\lmb_{ij}(1-g_ig_j),\label{eq:xixjlij} \\
&&x_i^{m_i} = \lmb_i(1-g_i^{m_i}).
\end{eqnarray}
It is well known that the set of monomials $\{x_1^{r_1}\cdots x_n^{r_n}g:0\le r_i<m_i,\;g\in G\}$ is a
basis of $H$. It is a Hopf algebra with comultiplication defined by
\begin{eqnarray}
	&& \Delta(g) = g\otimes g,\text{ for all } g\in G, \\
	&& \Delta(x_i) = g_i\otimes x_i+x_i\otimes 1.
\end{eqnarray}
The counit $\cou$ satisfies $\cou(g) = 1$, for all $g\in G$, and $\cou(x_i) = 0$.
Moreover, the antipode $\Ss$ is given by $\Ss(g) = g^{-1}$, for all $g\in G$, and
$\Ss(x_i)=-g_i^{-1}x_i$. We note that $\Ss^2(g) = g$ and $\Ss^2(x_i) = q_{ii}^{-1}x_i$.

Let $\mathbb{S}_n$ be the symmetric group on $n$ elements. For $\sigma \in\mathbb{S}_n$ let
$$
t_{\sigma} = x_{\sigma_{1}}^{m_{\sigma_{1}}-1}\cdots x_{\sigma_{n}}^{m_{\sigma_{n}}-1}
\sum_{g\in G}g.
$$
Note that $t_{\sigma}\ne 0$.

\begin{lemma}\label{lm:lqls} The following holds:
\begin{enumerate}
\item $\lmb_{ji}g_ig_j$ lies in the center of $H(\bfg,\bfch,\bflm)$ for $i\ne j$.  \label{lm:gigjz}
\item $\lmb_ig_i^{m_i}$ lies in the center of $H(\bfg,\bfch,\bflm)$.\label{lm:gimi}
\item $t_{\sigma}g = t_{\sigma}$, for all $g\in G$.
\item $t_{\sigma}x_{\sigma_{n}} = 0$.
\end{enumerate}
\end{lemma}

\begin{proof} (1) It is sufficient to see that $\lmb_{ji}g_ig_j$ commutes with $x_l$. If
$\lmb_{ji} = 0$ the result is clear. Assume that $\lmb_{ji}\neq 0$. Then, $\chi_i =
\chi_j^{-1}$, and thus
\begin{align*}
\lmb_{ji}g_ig_jx_l &= \lmb_{ji}\chi_l(g_i)\chi_l(g_j)x_lg_ig_j\\
&= \lmb_{ji}\chi_i(g_l)^{-1}\chi_j(g_l)^{-1}x_lg_ig_j\\
&= \lmb_{ji}x_lg_ig_j.
\end{align*}

\noindent (2) It is similar to \eqref{lm:gigjz}.

\noindent (3) It is immediate.

\noindent (4) We have:
\begin{align*}
t_{\sigma}x_{\sigma_{n}} &= x_{\sigma_{1}}^{m_{\sigma_{1}}-1}\cdots x_{\sigma_{n}}^{m_{\sigma_{n}} -1}
	\sum_{g\in G}gx_{\sigma_{n}}\\
& = x_{\sigma_{1}}^{m_{\sigma_{1}}-1}\cdots x_{\sigma_{n}}^{m_{\sigma_{n}}-1} x_{\sigma_{n}}
	\sum_{g\in G}\chi_{\sigma_{n}}(g)g \\
&= \lmb_{\sigma_{n}}x_{\sigma_{1}}^{m_{\sigma_{1}}-1}\cdots x_{\sigma_{n-1}}^{m_{\sigma_{n-1}}-1}
	(1-g_{\sigma_{n}}^{m_{\sigma_{n}}})\sum_{g\in G}\chi_{\sigma_{n}}(g)g,
\end{align*}
and the result follows by noticing that
\begin{align*}
(1-g_{\sigma_{n}}^{m_{\sigma_{n}}}) \sum_{g\in G}\chi_{\sigma_{n}}(g)g & = \sum_{g\in G}
\chi_{\sigma_{n}}(g) g - \sum_{g\in G} g_{\sigma_{n}}^{m_{\sigma_{n}}} \chi_{\sigma_{n}}(g)g\\
&= \sum_{g\in G}\chi_{\sigma_{n}}(g_{\sigma_{n}}^{m_{\sigma_{n}}}g) g_{\sigma_{n}}^{m_{
\sigma_{n}}} g-\sum_{g\in G}\chi_{\sigma_{n}}(g)g_{\sigma_{n}}^{m_{\sigma_{n}}}g = 0,
\end{align*}
since $\chi_{\sigma_{n}}(g_{\sigma_{n}}^{m_{\sigma_{n}}}) = q_{{\sigma_{n}} {\sigma_{n}}}^{
m_{\sigma_{n}}} = 1$.
\end{proof}

\begin{propo}\label{pr:pr4.2} $t_{\sigma}$ is a right integral.
\end{propo}

\begin{proof}
Let $M = (m_1-1)+\cdots+(m_n-1)$. Let
$$
\mathcal{A} = \{f\;:\;\{1,\ldots,M\}\to\{1,\ldots,n\}:\#f^{-1}(i) = m_i-1 \text{ for all } i\}.
$$
For $f\in\mathcal{A}$, let $x_f = x_{f(1)}x_{f(2)}\cdots x_{f(M)}$. We claim that if $f,h\in
\mathcal{A}$, then $x_f\sum_{g\in G}g = \beta x_h\sum_{g\in G}g$ for some $\beta\in\bk^*$. To
prove this claim, it is sufficient to check it when $f$ and $h$ differ only in $i,i+1$ for
some $1\le i<M$, that is, when $h = f\circ\tau_i$, where $\tau_i\in\Ss_M$ is the elementary
transposition $(i,i+1)$. But, in this case, we have:
\begin{align*}
x_f\sum_{g\in G}g &= x_{h\circ\tau_i}\sum_{g\in G}g \\
&= q_{h(i+1)h(i)}x_h\sum_{g\in G}g + \lambda_{h(i+1)h(i)}x_{h,\widehat i,\widehat{i+1}}
(1-g_{h(i)}g_{h(i+1)})\sum_{g\in G}g \\
&= q_{h(i+1)h(i)}x_h\sum_{g\in G}g
\end{align*}
where $x_{h,\widehat i,\widehat{i+1}} = x_{h(1)}\cdots x_{h(i-1)}x_{h(i+2)}\cdots x_{h(M)}$.
The second equality follows from relation~\eqref{eq:xixjlij} and item~\eqref{lm:gigjz} in the
previous Lemma. The Proposition follows now using items~(3) and (4) in the Lemma.
\end{proof}

Now we see that
\begin{itemize}
\item $\alpha(x_i) = 0$ (using Proposition~\ref{pr:pr4.2} and item~(2) of Lemma~\ref{lm:lqls}),
\item $\alpha(g) = \chi_1^{m_1-1}(g)\cdots\chi_n^{m_n-1}(g)$.
\end{itemize}
In particular, $\alpha(g_i) = q_{i1}^{m_1-1}\cdots q_{in}^{m_n-1}$. Since $\rho(h) =
\alpha(\Ss(\com h1))\Ss^2(\com h2)$, we have:
\begin{itemize}
\item $\rho(g) = \alpha(g^{-1})g$,
\item $\rho(x_i) = \alpha(g_i^{-1})q_{ii}^{-1}x_i =\displaystyle{\prod_{\substack{
{1\le j\le n}\\{j\neq i}}}q_{ij}^{1-m_j}x_i}$
\end{itemize}
Thus, as $\rho$ is an algebra map,
\begin{eqnarray*}
\rho(x_1^{r_1}\cdots x_n^{r_n}g)
   &=& q_{11}^{-r_1}\cdots q_{nn}^{-r_n}\alpha(g_1^{-r_1}\cdots g_n^{-r_n}g^{-1})
       \;x_1^{r_1}\cdots x_n^{r_n}g \\
   &=& \prod_{1\le i<j\le n}q_{ij}^{(1-m_j)r_i-(1-m_i)r_j}\alpha(g^{-1})
       \;x_1^{r_1}\cdots x_n^{r_n}g.
\end{eqnarray*}
So, the basis $\{x_1^{j_1}\cdots x_n^{j_n}g\}$ is made up of eigenvectors of $\rho$. Consider
the groups
$$
\bk^*\supseteq L_1 = \langle q_{11},\ldots,q_{nn},\alpha(G)\rangle\supseteq L_2 = \alpha(G).
$$
Using that $\rho(x_ig_i^{-1}) = q_{ii}^{-1}x_ig_i^{-1}$ and $\rho(g) = \alpha(g^{-1})g$, it is
easy to see that $L_1$ is the set of eigenvalues of $\rho$ and that the order of $\rho$ is the
l.c.m. of the numbers $m_1,\ldots,m_n$ and the order of the character $\alpha_{\vert G}\in\hat G$
(in particular, $\bk$ has a primitive $\ord_{\rho}$-th root of unity). As before, we decompose
$H = H(\bfg,\bfch,\bflm)$ as
\[
H = \bigoplus_{\omega\in L_1} H_{\omega},\quad\text{where } H_{\omega} = \{h\in H\;:\;\rho(h) =
\omega h\}.
\]
The following result is the version of Proposition~\ref{pr:gralc1} for the present context.
\begin{theor}\label{pr2:c1}
The following are equivalent:
\begin{enumerate}
\item\label{pr2:it1} $\bigoplus_{\omega\in L_1} H_{\omega}$ is strongly graded,
\item\label{pr2:it2} Each component $H_{\omega}$ contains an element in $G$,
\item\label{pr2:it3} $L_1 = L_2$,
\item\label{pr2:it4} $H$ is a crossed product $H_1\ltimes \bk L_1$.
\end{enumerate}
\end{theor}

\begin{proof} Clearly \eqref{pr2:it2} and \eqref{pr2:it3} are equivalent and
\eqref{pr2:it2}$\so$\eqref{pr2:it4}$\so$\eqref{pr2:it1}. Next we prove that
\eqref{pr2:it1}$\so$\eqref{pr2:it2}. Let $\omega\in L_1$. By Lemma~\ref{lm:nsg}, we know that
$\cou(H_{\omega})\neq 0$. Since $H_{\omega}$ has a basis consisting of monomials
$x_1^{r_1}\cdots x_n^{r_n}g$ and $\cou(x_i) = 0$, there must be an element $g\in G$ inside
$H_{\omega}$.
\end{proof}

\begin{remar}\label{rm:eqc}
	We next show that for liftings of Quantum Linear Spaces, the components in the decomposition
	$H = \bigoplus_{\omega\in L_1} H_{\omega}$ are equidimensional. In fact, in this case we can take the
	basis of $H$ given by
	$$
		\{(x_1g_1^{-1})^{r_1}\cdots(x_ng_n^{-1})^{r_n}g\;:\;0\le r_i<m_i,\;g\in G\}.
	$$
	Since $\rho(x_ig_i^{-1})=q_{ii}^{-1}x_ig_i^{-1}$, the map
	$$
		\theta\colon \bZ_{m_1}\times\cdots\times\bZ_{m_n}\times G\to\bk^*,
	$$
	taking $(r_1,\ldots,r_n,g)$ to the eigenvalue of $(x_1g_i^{-1})^{r_1}\cdots (x_ng_n^{-1})^{r_n}g$ with
	respect to $\rho$, is a well defined group homomorphism. From this it follows immediately that all the
	eigenspaces of $\rho$ are equidimensional.
\end{remar}

\section{Computing $H_1$}
Assume we are in the setting of the liftings of QLS. Suppose $H$ is a crossed product or, equivalently, that
$L_1 = L_2$. Then, there exist elements \mbox{$\gamma_1,\ldots,\gamma_n\in G$}, such that $\alpha(\gamma_i) =
q_{ii}$.  Set $\tilde\gamma_i = g_i^{-1}\gamma_i^{-1}$ and let $y_i = x_i\tilde\gamma_i$. It is immediate that
$y_i\in H_1$.  Let $N = \ker(\alpha_{\vert G}) \subseteq G$. It is easy to see that $H_1$ has a basis given by
$\{y_1^{r_1}\cdots y_n^{r_n}g:g\in N\}$. Furthermore, $H_1$ can be presented by generators $N,y_1,\ldots,y_n$
and relations
\begin{itemize}
\item $gy_i = \chi_i(g)y_ig$
\item $y_iy_j = q_{ij}\chi_j(\tilde\gamma_i)\chi_i^{-1}(\tilde\gamma_j)y_jy_i + \chi_j(
\tilde\gamma_i)\lambda_{ij}(\tilde\gamma_i\tilde\gamma_j-\gamma_i^{-1}\gamma_j^{-1})$
\item $y_i^{m_i} = \lambda_{i}\chi_i^{\frac{m_i(m_i-1)}2}(\tilde\gamma_i)(\tilde\gamma_i^{m_i}
-\gamma_i^{-m_i})$.
\end{itemize}
Notice that if $\lambda_{i}\neq 0$, then $\chi_i^{\frac{m_i(m_i-1)}2}(\tilde\gamma_i) = \pm
1$. We claim that
$$
	\lambda_{ij}\tilde\gamma_i\tilde\gamma_j,\quad \lambda_{ij}\gamma_i\gamma_j, \quad
	\lambda_i\tilde\gamma^{m_i}\quad\text{and}\quad \lambda_i\gamma^{m_i}
$$
belong to $kN$. It is clear that $\gamma^{m_i}\in N$, since $\alpha(\gamma^{m_i}) = q_{ii}^{m_i} = 1$. We now
prove the remaining part of the claim. Assume that $\lambda_{ij} \neq 0$. Then $\chi_i\chi_j = \varepsilon$.
Hence,

\begin{itemize}
	\item If $l\neq i,j$, then $\chi_l(g_ig_j) = q_{il}q_{jl} = q_{li}^{-1}q_{lj}^{-1} = \chi_i
		\chi_j(g_l^{-1}) = 1$.
	\item $q_{ii} = \chi_i(g_i) = \chi_j(g_i^{-1}) = q_{ij}^{-1} = q_{ji} = q_{jj}^{-1}$.
\end{itemize}

\noindent Thus, $m_i = \ord(q_{ii}) = \ord(q_{jj}) = m_j$, and then
$$
\chi_i^{m_i-1}(g_ig_j) \chi_j^{m_j-1}(g_ig_j) = (q_{ii}q_{ij}q_{ji}q_{jj})^{m_i-1} = 1
\quad\text{and}\quad \alpha(\gamma_i\gamma_j) = q_{ii}q_{jj} = 1.
$$
It is now immediate that $\alpha(g_ig_j) = \chi_1^{m_1-1}(g_ig_j)\cdots \chi_n^{m_n-1}(g_ig_j)
= 1$, and so
$$
	\alpha(\tilde\gamma_i\tilde\gamma_j) = \alpha(g_i^{-1}\gamma_i^{-1}g_j^{-1}\gamma_j^{-1}) =
	\alpha(g_jg_i)^{-1}\alpha(\gamma_j\gamma_i)^{-1} = 1.
$$
It remains to check that $\lambda_i\tilde\gamma^{m_i}\in kN$. Assume now that $\lambda_{i}\neq 0$. Then
$\chi_i^{m_i} = \varepsilon$. Thus,

\begin{itemize}
\item If $l\neq i$, then $\chi_l^{m_l-1}(g_i^{m_i}) = q_{il}^{(m_l-1)m_i} = q_{li}^{(1-m_l)m_i}
= \chi_i^{m_i}(g_l^{1-m_l}) = 1$.
\end{itemize}

\noindent Since $\chi_i^{m_i-1}(g_i^{m_i}) = q_{ii}^{m_i(m_i-1)} = 1 $, this implies that
$$
\alpha(g_i^{m_i}) = \chi_1^{m_1-1}(g_i^{m_i})\cdots \chi_n^{m_n-1}(g_i^{m_i}) = 1,
$$
and so
$$
\alpha(\tilde\gamma_i^{m_i}) = \alpha(\gamma_ig_i)^{-1} = \alpha(\gamma_i)^{-1}
\alpha(g_i)^{-1} = 1.
$$


\begin{thebibliography}{CD}
\bibitem[AG]{ag} N. Andruskiewitsch\ and\ M. Gra{\~n}a,
    \emph{Braided Hopf algebras over non-abelian finite groups},
    Bol. Acad. Nac. Cienc. (C\'ordoba), \textbf{63}, (1999), 45--78.
\bibitem[AS]{as} N. Andruskiewitsch\ and\ H.-J. Schneider,
    \emph{Pointed Hopf algebras},
    in New directions in Hopf algebras, Math. Sci. Res. Inst. Publ., \textbf{43}, 1--68,
    Cambridge Univ. Press, Cambridge, 2002.
\bibitem[A]{galgs} A. Marcus,
    \emph{Representation theory of group graded algebras},
    Nova Science Publishers Inc., Commack, NY, 1999.
\bibitem[GG]{gg} J.A. Guccione\ and\ J.J. Guccione,
    \emph{Hochschild cohomology of Frobenius algebras},
    Proc. Amer. Math. Soc., \textbf{132}, (2004), 5, 1241--1250 (electronic).
\bibitem[N]{n} W.D. Nichols,
    \emph{Bialgebras of type one},
    Comm. in Alg. {\bf 6} (1978), 1521--1552.
\bibitem[S]{S} H. J. Schneider,
    \emph{Lectures on Hopf Algebras},
    (1994)
\bibitem[R]{R} D.E. Radford,
    \emph{The order of the antipode of a finite dimensional Hopf algebra is finite},
    Amer. J. Math. {\bf 98} (1976), no.~2, 333--355.
\bibitem[RS]{rs} D. E. Radford\ and\ H.-J. Schneider,
    \emph{On the even powers of the antipode of a finite-dimensional Hopf algebra},
    J. Algebra {\bf 251} (2002), no.~1, 185--212.
\end{thebibliography}
\end{document}